\documentclass[12pt]{article}

\usepackage{amsmath,amsfonts,amsthm,enumerate,fullpage,young}
\usepackage{bbold} 
\usepackage{bm}

\def\S {\ensuremath{\mathcal{S}}}

\def\be {\begin{equation}}
\def\ee {\end{equation}}

\def\sgn		{\ensuremath{\mathrm{sgn}\,}}
\def\2			{\ensuremath{\mathbf{2}}}
\def\1			{\ensuremath{\mathbf{1}}}
\def\4			{\ensuremath{\mathbf{4}}}
\def\3			{\ensuremath{\mathbf{3}}}

\usepackage[vcentermath]{youngtab} 

\title{Young's Natural Representations of $\S_4$}
\author{Quinton Westrich}
\date{December 2008}

\begin{document}

\maketitle

\begin{abstract}
  We calculate all inequivalent irreducible representations of $\S_4$ by specifying the matrices for adjacent transpositions and indicating how to obtain general permutations in $\S_4$ from these transpositions. We employ standard Young tableaux methods as found in Sagan's \emph{The Symmetric Group: Representations, Combinatorial Algorithms, and Symmetric Functions} (2001).
\end{abstract}

\tableofcontents

\section{Introduction}

In this paper, we compute matrices of the adjacent transpositions $(1\,2),(2\,3),(3\,4)$ in $\S_4$ for all of Young's natural representations. These representations are irreducible and any other irreducible representation is equivalent to one of Young's natural representations. Moreover, the adjacent transpositions generate all of $\S_4$ so that all the matrices of Young's natural representations are easily obtained by merely multiplying the appropriate generators.

There are two key ideas employed in these calculations. The first is that partitions of $4$ correspond to irreducible representations of $\S_4$. Each partition is represented graphically by a Ferrer's diagram.
$$ \lambda_1 = (1,1,1,1) = \yng(1,1,1,1) \qquad \lambda_2 = (2,1,1) = \yng(2,1,1) \qquad \lambda_3 = (2,2) = \yng(2,2) $$
$$ \lambda_4 = (3,1) = \yng(3,1) \qquad \lambda_5 = (4) = \yng(4) $$
So there are five (inequivalent) irreducible representations of $\S_4$.

The second key idea is that given a partition of $4$, and hence a shape $\lambda$, the dimension of representation module corresponding to the shape $\lambda$, or the degree of the corresponding representation matrices, is equal to the number of standard fillings of the shape $\lambda$, denoted $f^\lambda$. This number is provided by the hook formula
$$ f^\lambda = \displaystyle\frac{n!}{\displaystyle\prod_{(i,j)\in\lambda}h_{ij}}. $$

Throughout this paper, we use the following convention in determining the representation of $(i\,i+1)e_t$ in terms of the basis of standard $\lambda$-polytabloids.
\begin{description}
	\item[Case 1.] $i$ and $i+1$ are in the same column of $t$. Then $$ (i,i+1) e_t = -e_t.$$
	\item[Case 2.] $i$ and $i+1$ are in the same row of $t$. Then we apply the Garnir element to obtain $(i,i+1)e_t$ as a linear combination of standard $\lambda$-polytabloids.
	\item[Case 3.] $i$ and $i+1$ are not in the same row or column of $t$. Then $(i,i+1)e_t$ is a standard $\lambda$-polytabloid.
\end{description}

\section{The Group $\S_4$}

There are $4!=24$ elements in $\S_4$. These elements fall into five conjugacy classes, each corresponding to a partition of $4$.
The conjugacy classes of $\S_4$, each labeled by its corresponding partition of $4$, are given below. \\
\begin{minipage}{0.8\textwidth}
	\begin{align*}
	  K_{(1,1,1,1)} &= \{\epsilon\} \\
	  K_{(2,1,1)} &= \{ (1\,2),(1\,3),(1\,4),(2\,3),(2\,4),(3\,4) \} \\
	  K_{(2,2)} &= \{ (1\,2)(3\,4),(1\,3)(2\,4),(1\,4)(2\,3) \} \\
	  K_{(3,1)} &= \{ (1\,2\,3),(1\,3\,2),(1\,2\,4),(1\,4\,2),(1\,3\,4),(1\,4\,3),(2\,3\,4),(2\,4\,3) \} \\
	  K_{(4)} &= \{ (1\,2\,3\,4),(1\,4\,3\,2),(1\,4\,2\,3),(1\,3\,2\,4),(1\,3\,4\,2),(1\,2\,4\,3) \}
	\end{align*}
\end{minipage}
\begin{minipage}{0.2\textwidth}
	\begin{align*}
	  |K_{(1,1,1,1)}| &= 1 \\
	  |K_{(2,1,1)}| &= 6 \\
	  |K_{(2,2)}| &= 3 \\
	  |K_{(3,1)}| &= 8 \\
	  |K_{(4)}| &= 6
	\end{align*}	
\end{minipage}

\vspace{0.3cm}

The elements of $\S_4$ are decomposed into adjacent transpositions below.
$$ \begin{array}{ll}
	\epsilon & (1\,2) \\
	(2\,3) & (3\,4) \\
	(1\,3) = (1\,2)(2\,3)(1\,2) & (1\,4) = (1\,2)(2\,3)(3\,4)(2\,3)(1\,2) \\
	(2\,4)=(2\,3)(3\,4)(2\,3) & (1\,2)(3\,4) \\
	(1\,3)(2\,4) = (1\,2)(2\,3)(1\,2)(2\,3)(3\,4)(2\,3) & (1\,4)(2\,3)=(1\,2)(2\,3)(3\,4)(2\,3)(1\,2)(2\,3) \\
	(1\,2\,3) = (1\,2)(2\,3) & (1\,3\,2) = (2\,3)(1\,2) \\
	(1\,2\,4) = (1\,2)(2\,3)(3\,4)(2\,3) & (1\,4\,2) = (2\,3)(3\,4)(2\,3)(1\,2) \\
	(1\,3\,4) = (1\,2)(2\,3)(3\,4)(1\,2) & (1\,4\,3) = (1\,2)(3\,4)(2\,3)(1\,2) \\
	(2\,3\,4) = (2\,3)(3\,4) & (2\,4\,3) = (3\,4)(2\,3) \\
	(1\,2\,3\,4) = (1\,2)(2\,3)(3\,4) & (1\,4\,3\,2) = (3\,4)(2\,3)(1\,2) \\
	(1\,4\,2\,3) = (1\,2)(2\,3)(1\,2)(2\,3)(3\,4)(2\,3)(1\,2) & (1\,3\,4\,2) = (2\,3)(3\,4)(1\,3) \\
	(1\,3\,2\,4) = (1\,2)(2\,3)(3\,4)(2\,3)(1\,2)(2\,3)(1\,2) & (1\,2\,4\,3) = (1\,2)(3\,4)(2\,3)
\end{array}$$ 
The homomorphism property of the representations then allows one to obtain the matrix for any element of $\S_4$ by simply multiplying the appropriate matrices according to the above computations.


\section{The $\S_4$-Module $\S^{\lambda_1}$}

\subsection{Standard $\lambda_1$-Tableaux}
Recall that $$ \lambda_1 = (1,1,1,1) = \yng(1,1,1,1). $$
So $$ f^{\lambda_1} = \displaystyle\frac{n!}{\displaystyle\prod_{(i,j)\in\lambda}h_{ij}}
 = \frac{24}{4\cdot 3 \cdot 2 \cdot 1} = 1. $$
Therefore, $\dim \S^{\lambda_1} = 1$. The standard filling is
$$ t = \young(1,2,3,4).$$
Denoting the representation by $X_{\lambda_1}$ we have immediately
$$ X_{\lambda_1}(\epsilon) = [1].$$

\subsection{Matrix for $(1\,2)$}

We have $(1\,2)e_t = e_{(1\,2)t} = -e_t $ by Case 1 and so
$$ X_{\lambda_1}(1\,2) = [-1]. $$

\subsection{Matrix for $(2\,3)$}

We have $(2\,3)e_t = e_{(2\,3)t} = -e_t $ by Case 1 and so
$$ X_{\lambda_1}(2\,3) = [-1]. $$

\subsection{Matrix for $(3\,4)$}

We have $(3\,4)e_t = e_{(3\,4)t} = -e_t $ by Case 1 and so
$$ X_{\lambda_1}(3\,4) = [-1]. $$


\section{The $\S_4$-Module $\S^{\lambda_2}$}
\subsection{Standard $\lambda_2$-Tableaux}

Recall that $$ \lambda_2 = (2,1,1) = \yng(2,1,1). $$
So $$ f^{\lambda_2} = \displaystyle\frac{n!}{\displaystyle\prod_{(i,j)\in\lambda}h_{ij}}
 = \frac{24}{4\cdot 1 \cdot 2 \cdot 1} = 3. $$
Therefore, $\dim \S^{\lambda_2} = 3$. The standard fillings are
$$ t_1 = \young(12,3,4) \qquad t_2 = \young(13,2,4) \qquad t_3 = \young(14,2,3).$$
Denoting the representation by $X_{\lambda_2}$ we have immediately
$$ X_{\lambda_2}(\epsilon) = \begin{bmatrix} 1 & 0 & 0 \\ 0 & 1 & 0 \\ 0 & 0 & 1\end{bmatrix}.$$

\subsection{Matrix for $(1\,2)$}

We have $$ (1\,2)t_1 = (1\,2)\young(12,3,4) $$
which is Case 2. We set $A=\{2,3,4\}$ and $B=\{1\}$. 
$$ \begin{array}{ccccc}
  (A',B'): & (2\,3\,4,1) & (1\,3\,4,2) & (2\,1\,4,3) & (2\,3\,1,4) \\
  & \downarrow\scriptstyle{\epsilon} & \downarrow\scriptstyle{(1\,2)} & \downarrow\scriptstyle{(1\,3)} & \downarrow\scriptstyle{(1\,4)} \\
  & \young(21,3,4) & \young(12,3,4) & \young(23,1,4) & \young(24,3,1) \\
  &&& \downarrow\scriptstyle{(1\,2)} & \downarrow\scriptstyle{(1\,2\,3)} \\
  &&& \young(13,2,4) & \young(14,2,3)
\end{array} $$
So $g_{A,B} = \epsilon - (1\,2) + (1\,2)(1\,3) = (1\,2\,3)(1\,4)$ and therefore
$$ (1\,2)e_{t_1} = e_{(1\,2)t_1} = e_{t_1} - e_{t_2} + e_{t_3}.$$

For the other two polytabloids we have Case 1.
$$ (1\,2)t_2 = (1\,2)\young(13,2,4) \qquad \Rightarrow \qquad (1\,2)e_{t_2} = -e_{t_2} $$
$$ (1\,2)t_3 = (1\,2)\young(14,2,3) \qquad \Rightarrow \qquad (1\,2)e_{t_3} = -e_{t_3} $$
Thus, the matrix for $(1\,2)$ is 
$$ X_{\lambda_2}(1\,2) = \begin{bmatrix} 1 & 0 & 0 \\ -1 & -1 & 0 \\ 1 & 0 & -1 \end{bmatrix}. $$

\subsection{Matrix for $(2\,3)$}

We have $$ (2\,3)t_1  = (2\,3)\young(12,3,4) = \young(13,2,4) = e_{t_2}$$
which is Case 3. Similarly,
$$ (2\,3)t_2 = (2\,3)\young(13,2,4) = \young(12,3,4) = e_{t_1}$$
which is Case 3. Finally,
$$ (2\,3)t_3 = (2\,3)\young(14,2,3) $$
which is Case 1. So $$ (2\,3)e_{t_3} = -e_{t_3}.$$
Therefore, the matrix for $(2\,3)$ is
$$ X_{\lambda_2}(2\,3) = \begin{bmatrix} 0 & 1 & 0 \\ 1 & 0 & 0 \\ 0 & 0 & -1 \end{bmatrix}. $$

\subsection{Matrix for $(3\,4)$}

We have $$ (3\,4)t_1 = (3\,4)\young(12,3,4) $$
which is Case 1. So $$ (3\,4)e_{t_1} = -e_{t_1}.$$
Next, $$ (3\,4)t_2 = (3\,4)\young(13,2,4) = \young(14,2,3) = t_3 ,$$
which is Case 3. So $$ (3\,4)e_{t_2} = e_{t_3}. $$ 
Finally, $$ (3\,4)t_3 = (3\,4)\young(14,2,3) = \young(13,2,4) = t_2, $$
which is Case 3. So $$ (3\,4)e_{t_3} = e_{t_2}. $$
Therefore, the matrix for $(3\,4)$ is 
$$ X_{\lambda_2}(3\,4) = \begin{bmatrix} -1 & 0 & 0 \\ 0 & 0 & 1 \\ 0 & 1 & 0 \end{bmatrix}. $$


\section{The $\S_4$-Module $\S^{\lambda_3}$}
\subsection{Standard $\lambda_3$-Tableaux}

Recall that $$ \lambda_3 = (2,2) = \yng(2,2). $$
So $$ f^{\lambda_3} = \displaystyle\frac{n!}{\displaystyle\prod_{(i,j)\in\lambda}h_{ij}}
 = \frac{24}{3\cdot 2 \cdot 2 \cdot 1} = 2. $$
Therefore, $\dim \S^{\lambda_3} = 2$. The standard fillings are
$$ t_1 = \young(12,34) \qquad t_2 = \young(13,24).$$
Denoting the representation by $X_{\lambda_3}$ we have immediately
$$ X_{\lambda_3}(\epsilon) = \begin{bmatrix} 1 & 0 \\ 0 & 1 \end{bmatrix}.$$

\subsection{Matrix for $(1\,2)$}

We have $$ (1\,2)t_1 = (1\,2)\young(12,34) $$
which is Case 2. Setting $A=\{2,3\}$ and $B=\{1\}$ we obtain
$$ \begin{array}{cccc}
	(A',B'): & (2\,3,1) & (1\,3, 2) & (1\,2,3) \\
	& \downarrow\scriptstyle{\epsilon} & \downarrow\scriptstyle{(1\,2)} & \downarrow\scriptstyle{(1\,3)} \\
	& \young(21,34) & \young(12,34) & \young(23,14) \\
	&&& \downarrow\scriptstyle{(1\,2)} \\
	&&& \young(13,24)
\end{array} $$
so that $e_{(1\,2)t_1} = e_{t_1} - e_{t_2}$.

We have $$(1\,2)t_2 = (1\,2)\young(13,24) $$
which is Case 1. So $(1\,2)e_{t_2} = -e_{t_2}$. Therefore, the matrix for $(1\,2)$ is 
$$ X_{\lambda_3}(1\,2) = \begin{bmatrix} 1 & 0 \\ -1 & -1 \end{bmatrix}. $$

\subsection{Matrix for $(2\,3)$}

We have Case 3 twice:
$$ (2\,3)t_1 = (2\,3)\young(12,34) = \young(13,24) = t_2$$
and 
$$ (2\,3)t_2 = (2\,3)\young(13,24) = \young(12,34) = t_1.$$
Therefore,
$$ (2\,3)e_{t_1} = e_{t_2} \qquad \mbox{ and } \qquad (2\,3)e_{t_2} = e_{t_1}.$$
Thus, the matrix for $(2\,3)$ is 
$$ X_{\lambda_3}(2\,3) = \begin{bmatrix} 0 & 1 \\ 1 & 0 \end{bmatrix}. $$

\subsection{Matrix for $(3\,4)$}

We have $$ (3\,4)t_1 = (3\,4)\young(12,34) = \young(12,43) $$
which is Case 2. Set $A=\{4\}$ and $B=\{2,3\}$. We obtain
$$ \begin{array}{cccc}
  (A',B'): & (4,2\,3) & (2,4\,3) & (3,2\,4) \\
  & \downarrow\scriptstyle{\epsilon} & \downarrow\scriptstyle{(2\,4)} & \downarrow\scriptstyle{(3\,4)} \\
  & \young(12,43) & \young(14,23) & \young(12,34) \\
  &&\downarrow\scriptstyle{(3\,4)} & \\
  && \young(13,24) &
\end{array} $$
which gives $$ e_{(3\,4)t_1} = -e_{t_2} + e_{t_1} .$$ 

Now $$(3\,4)t_2 = (3\,4) \young(13,24) $$
which is Case 1 so that $(3\,4)e_{t_2} = -e_{t_2}$. Thus, the matrix for $(3\,4)$ is 
$$ X_{\lambda_3}(3\,4) = \begin{bmatrix} 1 & 0 \\ -1 & -1 \end{bmatrix}.$$


\section{The $\S_4$-Module $\S^{\lambda_4}$}
\subsection{Standard $\lambda_4$-Tableaux}

Recall that $$ \lambda_4 = (3,1) = \yng(3,1). $$
So $$ f^{\lambda_4} = \displaystyle\frac{n!}{\displaystyle\prod_{(i,j)\in\lambda}h_{ij}}
 = \frac{24}{4\cdot 2 \cdot 1 \cdot 1} = 3. $$
Therefore, $\dim \S^{\lambda_4} = 3$. The standard fillings are
$$ t_1 = \young(134,2) \qquad t_2 = \young(124,3) \qquad t_3 = \young(123,4).$$
Denoting the representation by $X_{\lambda_4}$ we have immediately
$$ X_{\lambda_4}(\epsilon) = \begin{bmatrix} 1 & 0 & 0 \\ 0 & 1 & 0 \\ 0 & 0 & 1 \end{bmatrix}.$$

\subsection{Matrix for $(1\,2)$}

We have $$(1\,2)t_1 = (1\,2)\young(134,2)$$ which is Case 1 so that $(1\,2)e_{t_1} = -e_{t_1}$.

Now $$(1\,2)t_2 = (1\,2) \young(124,3) = \young(214,3)$$
which is Case 2. Set $A=\{2,3\}$ and $B=\{1\}$. Now
$$ \begin{array}{cccc}
  (A',B'): & (2\,3,1) & (1\,3,2) & (1\,2,3) \\
  & \downarrow\scriptstyle{\epsilon} & \downarrow\scriptstyle{(1\,2)} & \downarrow\scriptstyle{(1\,3)} \\
  & \young(214,3) & \young(124,3) & \young(234,1) \\
  &&& \downarrow\scriptstyle{(1\,2)} \\
  &&& \young(134,2)
\end{array} $$
which gives $e_{(1\,2)t_2} = e_{t_2} - e_{t_1}$.

Finally, $$(1\,2)t_3 = (1\,2) \young(123,4) = \young(213,4)$$
which is Case 2. Set $A=\{2,4\}$ and $B=\{1\}$. Then
$$ \begin{array}{cccc}
  (A',B'): & (2\,4,1) & (1\,4,2) & (1\,2,4) \\
  & \downarrow\scriptstyle{\epsilon} & \downarrow\scriptstyle{(1\,2)} & \downarrow\scriptstyle{(1\,4)} \\
  & \young(213,4) & \young(123,4) & \young(243,1) \\
  &&& \downarrow\scriptstyle{(1\,2)} \\
  &&& \young(143,2)
\end{array} $$
We have to find the Garnir element for the last tableau too. Set $\tilde A = \{4\}$ and $\tilde B=\{3\}$. Then
$$ \begin{array}{ccc}
  (\tilde{A'},\tilde{B'}): & (4,3) & (3,4) \\
  & \downarrow\scriptstyle{\epsilon} & \downarrow\scriptstyle{(3\,4)} \\
  & \young(143,2) & \young(134,2)
\end{array} $$
This provides
$$ e_{\scriptsize\young(143,2)} = e_{\scriptsize\young(134,2)} $$
as expected since the two tableaux are equivalent as tabloids. Altogether then 
$$ (1\,2) e_{t_3} = e_{t_3} - e_{t_1}.$$ 
This gives the matrix for $(1\,2)$ to be
$$ X_{\lambda_4}(1\,2) = \begin{bmatrix} -1 & -1 & -1 \\ 0 & 1 & 0 \\ 0 & 0 & 1 \end{bmatrix}. $$

\subsection{Matrix for $(2\,3)$}

We have $$(2\,3)t_1 = (2\,3)\young(134,2) = \young(124,3) = t_2 $$ as per Case 3. Likewise, 
$$ (2\,3)t_2 = (2\,3)\young(124,3) = \young(134,2) = t_1 $$ by Case 3. Thus,
$$ (2\,3)e_{t_1} = e_{t_2} \qquad \mbox{ and } \qquad (2\,3)e_{t_2} = e_{t_1}.$$
Lastly, $$(2\,3)t_3 = (2\,3)\young(123,4) $$
which calls for Case 2. Set $A=\{3\}$ and $B=\{2\}$, giving
$$ \begin{array}{ccc}
  (A',B'): & (3,2) & (2,3) \\
  & \downarrow\scriptstyle{\epsilon} & \downarrow\scriptstyle{(2\,3)} \\
  & \young(132,4) & \young(123,4)
\end{array} $$
so that $(2\,3)e_{t_3}=e_{t_3}$, again as expected when taking $t_3$ as a tabloid. Therefore
$$ X_{\lambda_4}(2\,3) = \begin{bmatrix} 0 & 1 & 0 \\ 1 & 0 & 0 \\ 0 & 0 & 1 \end{bmatrix}.$$

\subsection{Matrix for $(3\,4)$}

We have $$ (3\,4)t_1 = (3\,4) \young(134,2) = \young(143,2) $$
which is Case 2. Set $A=\{4\}$ and $B=\{3\}$. As for the cases above when $A$ and $B$ had one element, we obtain $e_{(3\,4)t_1} = e_{t_1}$.

The last two polytabloids obey Case 3. Thus,
$$ (3\,4) t_2 = (3\,4) \young(124,3) = \young(123,4) = t_3 \qquad\Rightarrow\qquad (3\,4)e_{t_2} = e_{t_3} $$
and
$$ (3\,4) t_3 = (3\,4) \young(123,4) = \young(124,3) = t_2 \qquad\Rightarrow\qquad (3\,4)e_{t_3} = e_{t_2} $$
so that the matrix for $(3\,4)$ is
$$ X_{\lambda_4}(3\,4) = \begin{bmatrix} 1 & 0 & 0 \\ 0 & 0 & 1 \\ 0 & 1 & 0 \end{bmatrix}. $$


\section{The $\S_4$-Module $\S^{\lambda_5}$}

\subsection{Standard $\lambda_5$-Tableaux}

Recall that $$ \lambda_5 = (4) = \yng(4). $$
So $$ f^{\lambda_5} = \displaystyle\frac{n!}{\displaystyle\prod_{(i,j)\in\lambda}h_{ij}}
 = \frac{24}{4\cdot 3 \cdot 2 \cdot 1} = 1. $$
Therefore, $\dim \S^{\lambda_5} = 1$. The standard filling is
$$ t = \young(1234).$$
Denoting the representation by $X_{\lambda_5}$ we have immediately
$$ X_{\lambda_5}(\epsilon) = [1].$$

\subsection{Matrices for $(1\,2)$, $(2\,3)$, and $(3\,4)$}

In every case, we obtain Case 2 with $A$ and $B$ being singletons. As above, we obtain
$$(1\,2)e_t = e_t \qquad (2\,3)e_t=e_t \qquad \mbox{and} \qquad (3\,4)e_t = e_t.$$
(Again, this can be easily seen by viewing $t$ as its corresponding tabloid $\{t\}$.) Therefore the matrices are
$$ X_{\lambda_5}(1\,2) = [1] \qquad X_{\lambda_5}(2\,3) = [1] \qquad \mbox{and} \qquad X_{\lambda_5}(3\,4) = [1].$$


\section{Conclusion}

We see that the representations corresponding to $\lambda_1=\scriptsize\yng(1,1,1,1)$ and $\lambda_2 = \scriptsize\yng(4)$ are the sign and trivial representations, respectively. 

We can compare the characters we obtain with those obtained in a previous work. There we had
	  \begin{center}
	  \begin{tabular}{r|ccccc}
	    $\S_4$ & $K_{(1,1,1,1)}$ & $K_{(2,1,1)}$ & $K_{(2,2)}$ & $K_{(3,1)}$ & $K_{(4)}$ \\
	    \hline\hline
	    trivial & $1$ & $1$ & $1$ & $1$ & $1$ \\
	    sign & $1$ & $-1$  & $1$ & $1$ & $-1$ \\
	    $\chi$ & $3$ & $1$ & $-1$ & $0$ & $-1$ \\
	    $\chi\hat\otimes\chi_{\sgn}$ & $3$ & $-1$ & $-1$ & $0$ & $1$ \\
	    $\psi$ & $2$ & $0$ & $2$  & $-1$ & $0$
	  \end{tabular}
	  \end{center}
To obtain the characters from Young's natural representations, we notice that each row in the table is completely characterized by its values in the first and second columns, which correspond to the character on the identity and on an adjacent transposition. Therefore we can conclude that 
\begin{align*}
	\chi &\leftrightarrow X_{\lambda_4} \leftrightarrow \yng(3,1) \\
	\chi\hat\otimes\chi_{\sgn} &\leftrightarrow X_{\lambda_2} \leftrightarrow \yng(2,1,1) \\
	\psi &\leftrightarrow X_{\lambda_3} \leftrightarrow \yng(2,2) \\
\end{align*}
Using the lexicographic order on columns then, we have a standard presentation of the character table of $\S_4$.
	  \begin{center}
	  \begin{tabular}{ccr|ccccc}
	    && $\S_4$ & $K_{(1,1,1,1)}$ & $K_{(2,1,1)}$ & $K_{(2,2)}$ & $K_{(3,1)}$ & $K_{(4)}$ \\
	    \hline\hline
	    \rule{0pt}{0.9cm} $\scriptsize\yng(1,1,1,1)$ & $X_{\lambda_1}$ & sign & $1$ & $-1$  & $1$ & $1$ & $-1$ \\[0.5cm]
	    $\scriptsize\yng(2,1,1)$ & $X_{\lambda_2}$ & $\chi\hat\otimes\chi_{\sgn}$ & $3$ & $-1$ & $-1$ & $0$ & $1$ \\[0.4cm]
	    $\scriptsize\yng(2,2)$ & $X_{\lambda_3}$ & $\psi$ & $2$ & $0$ & $2$  & $-1$ & $0$ \\[0.3cm]
	    $\scriptsize\yng(3,1)$ & $X_{\lambda_4}$ & $\chi$ & $3$ & $1$ & $-1$ & $0$ & $-1$ \\[0.3cm]
	    $\scriptsize\yng(4)$ & $X_{\lambda_5}$ & trivial & $1$ & $1$ & $1$ & $1$ & $1$ 
	  \end{tabular}
	  \end{center}

\end{document}